\documentclass[final,1p,times]{elsarticle}

\usepackage{graphicx}
\usepackage{color}
\newtheorem{thm}{Theorem}

\usepackage{amssymb,amsmath}

\usepackage{multirow}
\usepackage{stackrel}

\newcommand{\e}{\mathrm{e}}
\newcommand{\dd}{\mathrm{d}}

\newdefinition{rmk}{Remark}
\newtheorem{example}{Example}
\journal{XXX}
\usepackage[all]{xy}

\begin{document}

\begin{frontmatter}



\title{A probabilistic analysis of a Beverton-Holt type discrete model: Theoretical and computing analysis}


\author[IMM]{J.-C. Cort\'es\corref{cor1}}
\ead{jccortes@imm.upv.es}

\author[deusto1,deusto2]{A. Navarro-Quiles}
\ead{annaqui@doctor.upv.es}

\author[IMM]{J.-V. Romero}
\ead{jvromero@imm.upv.es}

\author[IMM]{M.-D. Rosell\'{o}}
\ead{drosello@imm.upv.es}

\address[IMM]{Instituto Universitario de Matem\'{a}tica Multidisciplinar,\\
Universitat Polit\`{e}cnica de Val\`{e}ncia,\\
Camino de Vera s/n, 46022, Valencia, Spain}

\address[deusto1]{DeustoTech, University of Deusto, 48007 Bilbao,\\
 Basque Country, Spain}
\address[deusto2]{ Facultad de Ingenieria, Universidad de Deusto, \\
Avda.Universidades, 24, 48007, Bilbao, Basque Country, Spain.}

\cortext[cor1]{Corresponding author. Phone number: +34--963877000 (ext.~88289)}

\begin{abstract}
In this paper a randomized version of the Beverton-Holt type discrete model is proposed. Its solution stochastic process and the random steady state are determined.  Its first probability density function and second probability density function are obtained by means of the random variable transformation method, providing a full probabilistic description of the solution. Finally, several numerical examples are shown.

\end{abstract}


\end{frontmatter}


\section{Introduction} \label{introduction}

One of the most interesting problems in population dynamics is  modelling  the changes  in the size of the population over the time \cite{molles2007ecology}. The breeding population of a species has been mathematically described using  both continuous and discrete models. In the former case, the corresponding models are mainly based on differential equations while difference equations are the most used mathematical tools  for the latter approach. 

Explicit first-order differential equations are formulated in the form
\begin{equation}\label{edo_general}
x'(t)=f(x(t),t),
\end{equation}
where the unknown, $x(t)$, represents  the number or percentage of individuals of the population under study, while the right-hand size term, $f(x(t),t)$, is a specific function depending on the proposed model. In the the  simplest model, usually termed Malthusian or exponential model, $f(x(t),t)= \alpha x(t)$. This model is satisfactory as long as the population is not too large with respect to its environment (resources). Otherwise, the accuracy of the model fails and, when $\alpha>0$, the Malthusian model predicts an exponential increase in the population as time goes on, which is unrealistic because of resources,  like food, are always limited.  This drawback is overcame by  the most popular continuous model of the growth of population, usually referred to as  the logistic o Verhulst-Pearl model formulated via a nonlinear differential equation where $f(x(t),t)=x(t)(\alpha-\beta x(t))$, being $\alpha,\beta>0$. In this model $\alpha$ represents the rate of the growth of the population if the resources were unlimited and the individuals did not affect  one another, and the term $-\beta x^2(t)$ can be interpreted as the negative effect on the growth due to crowdedness and limited resources. This model was first published by Pierre Verhulst \cite{ver,ver2}. This  equation has been extensively used to model different problems in a variety of  realms including Social Sciences, Biology,  Ecology, etc. \cite{KWASNICKI20135,Banasiak}, showing its great flexibility in applications. It is well-known that the general solution of logistic model is given by $x(t)=(\alpha / \beta)/(1+\left(\e^{-\alpha t}/(\gamma \beta)\right))$, where the parameter $\gamma$ is a real number that can be determined once an initial condition (the initial population) has been set.

Many times only data is available for discrete times, so it is interesting to have the solution for these time instants. In that case, the application of discrete models is natural and difference equations emerge. 
In this context, a very useful discrete equation in ecology is the non-autonomous Beverton-Holt equation\cite{Beverton}
\[
  x_{n+1} = \frac{\nu K_n x_n}{K_n+(\nu -1) x_n}, \qquad n \geq 0,
\]
where $\nu \geq 1$, $K_n>0$ and $x_0>0$. The  Beverton-Holt equation is of great importance in the fishery industry concerning the grow and exploitation of some species.
The sequence $K_n$ is the carrying capacity and $\nu$ is the inherent growth rate. Sometimes the carrying capacity can be considered constant\cite{Beverton}, it is $K_n=K$. This leads to the following difference equation
\begin{equation}\label{eqBH}
  x_{n+1}= \frac{a x_n}{1+ b x_n},
\end{equation} 
where $a=\nu \geq 1$ and $b=\frac{1}{K}(\nu -1)=\frac{1}{K}(a-1) \geq 0$.

An interesting property is that Beverton-Holt equation \eqref{eqBH} can be derived from logistic equation, as we will show below. 
Denoting an arbitrary  time instant by  $n$, and using the standard notation $x_n$ for the solution of the logistic model at  $n$, we clearly have  $x_n=(\alpha / \beta)/(1+ \e^{-\alpha n}/(\gamma \beta))$. Now, let us observe that
\[
x_{n+1}=\frac{\alpha / \beta}{1+\left(\e^{-\alpha(n+1)}/(\gamma \beta)\right)}
=\frac{\e^{\alpha}\alpha / \beta}{\e^{\alpha}+\left(\e^{-\alpha n}/(\gamma \beta)\right)}
=\frac{\e^{\alpha}\alpha / \beta}{1+\left(\e^{-\alpha n}/(\gamma \beta)\right)+\e^{\alpha}-1}.
\]
From this latter expression, we can derive a difference equation associated with $x_n$. Indeed, let us observe that
\[
\begin{array}{ccl}
x_{n+1} & =& \dfrac{\e^{\alpha}\alpha / \beta}{1+ \e^{-\alpha n}/(\gamma \beta)}
\dfrac{1+ \e^{-\alpha n}/(\gamma \beta)}{1+\left(\e^{-\alpha n}/(\gamma \beta)\right)+\e^{\alpha}-1}
\\
\\
 & = & \dfrac{\e^{\alpha}\alpha / \beta}{1+ \e^{-\alpha n}/(\gamma \beta)}
\dfrac{1}{1+\dfrac{\e^{\alpha} -1}{1+ \e^{-\alpha n}/( \gamma \beta)}}
\\
\\
 & = & \e^{\alpha} \dfrac{\alpha / \beta}{1+ \e^{-\alpha n}/(\gamma \beta)}
 \dfrac{1}{1+ \beta/\alpha \left( \e^{\textcolor{blue}{\alpha}}-1 \right) \dfrac{\alpha/\beta}{1+(\e^{-\alpha n})/(\gamma \beta )}}
 \\
\\
 & = & \e^{\alpha}
 x_n \dfrac{1}{1+ \beta/\alpha \left( \e^{\textcolor{blue}{\alpha}}-1 \right) x_n}
 \\
\\
 & = &  \dfrac{\e^{\alpha}
 x_n}{1+ \dfrac{\beta}{\alpha} \left( \e^{\textcolor{blue}{\alpha}}-1 \right) x_n}.
\end{array}
\]
Therefore,
\begin{equation} \label{probdet}
x_{n+1}=\dfrac{\e^{\alpha}
 x_n}{1+ \dfrac{\beta}{\alpha} \left( \e^{\textcolor{blue}{\alpha}}-1 \right) x_n}.
\end{equation}
Observe that this model has the structure of Eq. \eqref{eqBH}, where $a=\e^{\alpha}>1$ and $b=\frac{\beta}{\alpha}\left(\e^{\alpha}-1 \right)>0 $, that belongs to the so-called Beverton-Holt type (sometimes referred to as Pielou logistic equation \cite{De_La_Sen}).

Beverton-Holt type equation \eqref{eqBH} is a nonlinear difference equation. By letting  $x_n=1/z_n$ it is transformed into the linear equation
\begin{equation}
 z_{n+1} = \frac{1}{a} z_n + \frac{b}{a},
\end{equation}
whose solution, taking the initial condition $z_0=1/c$, is
\[ 
 z_n =
 \left\{
   \begin{array}{llll}
   \displaystyle 
     \left[ \frac{1}{c} - \frac{b}{a-1} \right] a^{-n} + \frac{b}{a-1}, &&& \text{if } a \neq 1, \\[0.6cm]
   \displaystyle 
      \frac{1}{c} + b n, &&& \text{if } a = 1.
   \end{array}
 \right.
\]
Or equivalently, in terms of unknown $x_n$ can be written as  
\begin{equation}  \label{soldet}
 x_n =
 \left\{
   \begin{array}{llll}
   \displaystyle 
     \frac{a^n (a -1)}{b a^n + \frac{1}{c} (a -1)-b}, &&& \text{if } a \neq 1, \\[0.6cm]
   \displaystyle 
      \frac{1}{\frac{1}{c}+b n} &&& \text{if } a = 1.
   \end{array}
 \right.
\end{equation}

An important magnitude is the steady state, which is given by
\begin{equation} \label{estacionario}
\lim_{n\rightarrow \infty}  x_n =
 \left\{
   \begin{array}{llll}
   \displaystyle 
     \frac{a -1}{b}, &&& \text{if } a \neq 1, \\[0.6cm]
   \displaystyle 
      0 &&& \text{if } a = 1.
   \end{array}
 \right.
\end{equation}

In real problems coefficients and initial conditions are not usually known exactly.  This may be due to measurement errors or the
inherent complexity associated to their own nature. So, it  seems more realistic to consider that parameters and initial conditions are random variables instead of deterministic magnitudes.

Denoting random variables by capital letters, $A=\e^{\alpha}>1$ and $B=\frac{\beta}{\alpha}\left(\e^{\alpha}-1 \right)>0 $, the  random Beverton-Holt model \eqref{probdet} can be written as
\begin{equation} \label{Pielourand}
\left\{
\begin{array}{l}
   \displaystyle X_{n+1}= \frac{A X_n}{1+B X_n}, \qquad n=0,1,2,\ldots , \\[0.4cm]
   X_{0}=C,
\end{array}
\right.
\end{equation}
where all the input parameters $A$,$B$ and $C$ are assumed to be absolutely continuous random variables defined on a common complete probability space 
($\Omega, \mathcal{F},\mathbb{P}$). Hereinafter, we assume that 
$\mathbb{P}\left[ \{ \omega \in \Omega:\, A(\omega)>1 \} \right]=1$,
$\mathbb{P}\left[ \{ \omega \in \Omega:\, B(\omega)>0 \} \right]=1$
and $\mathbb{P}\left[ \{ \omega \in \Omega:\, C(\omega)>0 \} \right]=1$. Moreover, for the sake of completeness, we will suppose that $A,B,C$ are dependent random variables whose joint PDF is $f_{C,A,B}(c,a,b)$.

Solving a random difference  equation means not only to calculate its exact solution  stochastic process but also to determine its main statistical functions like the mean and the variance Additionally, an important goal is the computation of the 1-PDF, say $f^X_1(x,n)$, since from it  we can compute all the one-dimensional statistical moments of any order. In particular, from the 1-PDF we can straightforwardly determine the mean function
\[
\mathbb{E}\left[X_n\right]=\int_{\mathbb{R}} x f^X_1(x,n)\mathrm{d}x,
\] 
and  the variance function
\[
\sigma^2\left[X_n\right]=\int_{\mathbb{R}} x^2 f^X_1(x,n)\mathrm{d}x-\left(\mathbb{E}\left[X_n\right]\right)^2.
\] 

The computation of the 1-PDF allows us to obtain a complete statistical description of the solution stochastic process at every time instant or period, say $n$. In general, a more challenging problem is to obtain the rest $n$-dimensional PDFs of the solution stochastic process for $n\geq 2$ (the so-called fidis: finite-dimensional distributions) because it usually involves complex computations. These functions provide us important probabilistic information. For example, the 2-PDF, $f^X_2(x_1,n_1;x_2,n_2)$, allows us to obtain a complete probabilistic description of the solution stochastic process at every arbitrary pair of periods, say $n_1$ and $n_2$. In particular, from the 2-PDF we can compute the correlation function. The correlation function gives a measure of linear statistical interdependence between $X_{n_1}$ and $X_{n_2}$
\begin{equation}\label{correlacion}
\Gamma_X(n_1,n_2)=\mathbb{E}\left[ X_{n_1} X_{n_2}  \right]=\int_{\mathbb{R}^2} x_1x_2 f^X_2(x_1,n_1;x_2,n_2)\,\mathrm{d}x_1 \mathrm{d}x_2.
\end{equation}
Furthermore, $\Gamma_X(n_1,n_2)$ together with the mean function, permits us the computation of the covariance function
\begin{equation}\label{covarianza}
C_X(n_1,n_2)=\Gamma_X(n_1,n_2)-\mathbb{E}\left[X_{n_1}\right]\mathbb{E}\left[X_{n_2}\right].
\end{equation}
The so-called RVT method is a powerful method  that has been recently used  to determine the PDF of the solution of some relevant differential and difference equations \cite{CORTES2018190,CASABAN2017396,CORTES2017225,CORTES2017150}.
The RVT technique permits to compute the PDF of a random vector which results from mapping another random vector whose PDF is  known. For the sake of completeness, now we state the multidimensional version of the RVT technique that will be extensively applied throughout this paper.

\begin{thm} (Multidimensional version, \cite[pp. 24--25]{soong}). \label{R.V.T.method2}
Let $\mathbf{U}=(U_1,\ldots,U_n)^{\top}$ and $\mathbf{V}=(V_1,\ldots,V_n)^{\top}$ be two $n$-dimensional absolutely continuous random vectors. Let $\mathbf{r}: \mathbb{R}^n \rightarrow \mathbb{R}^n$ be a one-to-one deterministic transformation of $\mathbf{U}$ into $\mathbf{V}$, i.e., $\mathbf{V}=\mathbf{r}(\mathbf{U})$. Assume that $\mathbf{r}$ is continuous in $\mathbf{U}$ and has continuous partial derivatives with respect to $\mathbf{U}$. Then, if $f_{\mathbf{U}}(\mathbf{u})$ denotes the joint PDF of vector $\mathbf{U}$, and $\mathbf{s}=\mathbf{r}^{-1}=(s_1(v_1,\ldots,v_n),\ldots,s_n(v_1,\ldots,v_n))^{\top}$ represents the inverse mapping of $\mathbf{r}=(r_1(u_1,\ldots,u_n),\ldots,r_n(u_1,\ldots,u_n))^{\top}$, the joint PDF of vector $\mathbf{V}$ is given by
\begin{equation}\label{RVTformula}
f_{\mathbf{V}}(\mathbf{v})=f_{\mathbf{U}}\left(\mathbf{s}(\mathbf{v})\right) \left| J_n \right|,
\end{equation}
where $\left| J_n \right|$ is the absolute value of the Jacobian, which is defined by
\begin{equation}\label{jacobiano_general}
J_n=\det \left( \frac{\partial \mathbf{s}^\top}{\partial \mathbf{v}}\right)
=
\det
\left(
\begin{array}{ccc}
\dfrac{\partial s_1(v_1,\ldots, v_n)}{\partial v_1} & \cdots & \dfrac{\partial s_n(v_1,\ldots, v_n)}{\partial v_1}\\
\vdots & \ddots & \vdots\\
\dfrac{\partial s_1(v_1,\ldots, v_n)}{\partial v_n} & \cdots & \dfrac{\partial s_n(v_1,\ldots, v_n)}{\partial v_n}\\
\end{array}
\right)
.
\end{equation}
\end{thm}

The layout of the paper is as follows. In Section \ref{sec_Beverton-Holt} we compute explicit expressions of the 1-PDF (Subsection\ref{subseccion_1PDF}) and the 2-PDF (Subsection\ref{subseccion_2PDF}) of the solution stochastic process of Beverton-Holt model in terms of the joint PDF of the random inputs. In Section \ref{sec_steady_state}, the PDF of the equilibrium of Beverton-Holt model is determined. Section \ref{sec_ejemplos} is addressed to illustrate of our theoretical results via some numerical experiments assuming a wide range of statistical distributions for the input model parameters. Conclusions are drawn in Section \ref{sec_Conclusiones}. 

\section{Solving the randomized Beverton-Holt model}\label{sec_Beverton-Holt}

The goal of this section is to obtain the two first PDFs (1-PDF and 2-PDF) of the solution of \eqref{Pielourand}, say $f_1^{X}(x,n)$ and $f_2^{X}(x_1,n_1;x_2,n_2)$, respectively. As $A$ is an absolutely continuous random variable, then  $\mathbb{P}\left[ \{ \omega \in \Omega:\, A(\omega)=1 \} \right]=0$, for all event $\omega\in \Omega$. As a consequence, taking into account that the solution of deterministic problem \eqref{probdet} (with $a=\e^{\alpha}>1$ and $b=\frac{\beta}{\alpha}\left(\e^{\alpha}-1 \right)>0 $) is given by \eqref{soldet}, we obtain that  the solution of random Beverton-Holt model \eqref{Pielourand} is given by
\begin{equation} \label{sol}
X_n
= \frac{ A^n (A-1)}{ B A^n+ \frac{1}{C} (A-1)-  B}
,\quad n=0,1,\ldots
\end{equation}
and this solution is well-defined.

\subsection{Computing the 1-PDF of the solution stochastic process}\label{subseccion_1PDF}

As previously indicated, in this section we will obtain the 1-PDF of \eqref{sol}  using the RVT method. To do this, given an arbitrary but fixed period $n$, we will apply Theorem~\ref{R.V.T.method2} for the following choice of mapping $\mathbf{r}$ 
\[
\begin{array}{ccccl}
y_1 & = & r_1\left( c,a,b \right) &=& \displaystyle \frac{a^n (a-1)}{ba^n+ \frac{1}{c} (a-1)-b},\\[0.2cm]
y_2 & = & r_2\left( c,a,b \right) &=& a,\\[0.2cm]
y_3 & = & r_3\left( c,a,b \right) &=& b.
\end{array}
\]
 So, the inverse mapping  $\mathbf{s}$ of  mapping $\mathbf{r}$ is given by
\[
\begin{array}{lllll}
c & = & s_1 \left(y_1,y_2,y_3  \right) &=& \displaystyle
\frac{y_1 \left(y_2 -1 \right)}{y_2^n \left( y_2 -1 \right) - y_3 y_1 \left( y_2^n -1 \right)},\\[0.2cm]
a & = & s_2 \left(y_1,y_2,y_3 \right) & =&y_2 ,\\[0.2cm]
b & = & s_3 \left(y_1,y_2,y_3 \right) &=& y_3, 
\end{array}
\]
and the absolute value of its Jacobian is
\[
\left | J_3 \right | = \left |  \frac{\partial s_1}{\partial y_1}\right|=
\left |  \frac{\left( y_2 -1 \right)^2 y_2^n}{\left( y_2^n \left( y_2 -1 \right) - y_3 y_1 \left( y_2^n -1 \right)  \right)^2}\right|,
\]
which is different from zero since $y_2=a=\e^{\alpha}>1$. Observe that $0\neq y_2\neq 1$  since $A$ is an absolutely continuous RV.

Applying  Theorem~\ref{R.V.T.method2}, the PDF of the random vector $(Y_1,Y_2,Y_3)$ defined by mapping $\mathbf{r}$ is
\begin{equation} \label{PDFconj}
\begin{array}{ccl}
f_{Y_1,Y_2,Y_3}\left( y_1,y_2,y_3\right)&=&
\displaystyle
f_{C,A,B}
\left(
\frac{y_1 \left(y_2 - 1 \right)}{y_2^n \left( y_2 -1 \right) - y_3 y_1 \left( y_2^n -1 \right)},
y_2,y_3
\right)\\[0.4cm]
& & \times  \displaystyle
\left |  \frac{\left( y_2 -1 \right)^2 y_2^n}{\left( y_2^n \left( y_2 -1 \right) - y_3 y_1 \left( y_2^n -1 \right)  \right)^2}\right|.
\end{array}
\end{equation}
Now, marginalizing expression \eqref{PDFconj} with respect to $A$ and $B$ and being $n$ arbitrary, we obtain the 1-PDF of $X_n$
\begin{equation}
\label{f1xn}
f_1^{X}(x,n)= \displaystyle
\iint_{\mathbb{R}^2} 
f_{C,A,B}
\left(
\frac{x \left(a -1 \right)}{a^n \left( a -1 \right) - b x \left( a^n -1 \right)},
a,b
\right)
\left |  \frac{\left( a -1 \right)^2 a^n}{\left( a^n \left( a -1 \right) - b x \left( a^n -1 \right)  \right)^2} \right| \dd a \, \dd b.
\end{equation}

\subsection{Computing the 2-PDF of the solution stochastic process}\label{subseccion_2PDF}

To compute the 2-PDF of the solution stochastic process $X_n$, an analogous  reasoning exhibited in the previous subsection will be applied. Let us consider two periods of time $n_1,n_2\geq 0$ fixed, being $n_1\neq n_2$, and then we apply Theorem~\ref{R.V.T.method2} choosing the mapping $\mathbf{r}$ as follows
\[
\begin{array}{ccccl}
y_1 & = & r_1\left( c,a,b \right) &=& \displaystyle \frac{a^{n_1} (a-1)}{ba^{n_1}+ \frac{1}{c} (a-1)-b},\\[0.5cm]
y_2 & = & r_2\left( c,a,b \right) &=&  a,\\[0.2cm]
y_3 & = & r_3\left( c,a,b \right) &=& \displaystyle \frac{a^{n_2} (a-1)}{ba^{n_2}+ \frac{1}{c} (a-1)-b}.
\end{array}
\]
The inverse mapping $\mathbf{s}$  of $\mathbf{r}$  is given by
\[
\begin{array}{lllll}
c & = & s_1 \left(y_1,y_2,y_3  \right) &=& \displaystyle
\frac{y_1 y_3 \left(y_2^{n_1} -y_2^{n_2} \right)}{y_1 y_2^{n_2} \left( -1+y_2^{n_1} \right) - y_3 y_2^{n_1} \left( -1 + y_2^{n_2} \right)},\\[0.5cm]
a & = & s_2 \left(y_1,y_2,y_3 \right) & =&y_2 ,\\[0.2cm]
b & = & s_3 \left(y_1,y_2,y_3 \right) &=& \displaystyle
\frac{(-1+y_2)\left(y_3 y_2^{n_1} -y_1 y_2^{n_2} \right)}{y_1 y_3 \left( y_2^{n_1}-y_2^{n_2} \right)}, 
\end{array}
\]
which absolute value of Jacobian
\[
\left | J_3 \right | = \left |  \frac{\partial s_1}{\partial y_1} \frac{\partial s_3}{\partial y_3}-\frac{\partial s_1}{\partial y_3}\frac{\partial s_3}{\partial y_1}\right|=
\left |  \frac{  (-1+y_2)y_2^{n_1+n_2}(y_2^{n_1}-y_2^{n_2})   }{(y_1y_2^{n_2}(-1+y_2^{n_1})-y_3 y_2^{n_1}(-1+y_2^{n_2}))^2}  \right|,
\]
is different from zero since $y_2>1$ and $0\neq y_2 \neq 1$, as it was justified previously. Applying the RVT technique, Theorem \ref{R.V.T.method2}
\[
\begin{array}{lcl}
f_{Y_1,Y_2,Y_3}(y_1,y_2,y_3)&=&\displaystyle f_{C,A,B}\left(\frac{y_1 y_3 \left(y_2^{n_1} -y_2^{n_2} \right)}{y_1 y_2^{n_2} \left( -1+y_2^{n_1} \right) - y_3 y_2^{n_1} \left( -1 + y_2^{n_2} \right)}, y_2, \frac{(-1+y_2)\left(y_3 y_2^{n_1} -y_1 y_2^{n_2} \right)}{y_1 y_3 \left( y_2^{n_1}-y_2^{n_2} \right)} \right)\\
\\
& &\displaystyle \times \left |  \frac{  (-1+y_2)y_2^{n_1+n_2}(y_2^{n_1}-y_2^{n_2})   }{\left( y_1y_2^{n_2}(-1+y_2^{n_1})-y_3 y_2^{n_1}(-1+y_2^{n_2})\right)^2}  \right|
\end{array}
\]
Now, marginalizing with respect to $A$ and taking $n_1$ and $n_2$ arbitrary, the 2-PDF is
\begin{equation}\label{2PDFconj}
\begin{array}{lcl}
f_2^X(x_1,n_1;x_2,n_2)&=& \displaystyle \int_{\mathbb{R}} f_{C,A,B}\left(\frac{x_1 x_2 \left(a^{n_1} -a^{n_2} \right)}{x_1 a^{n_2} \left( -1+a^{n_1} \right) - x_2 a^{n_1} \left( -1 + a^{n_2} \right)}, a, \frac{(-1+a)\left(x_2 a^{n_1} -x_1 a^{n_2} \right)}{x_1 x_2 \left( a^{n_1}-a^{n_2} \right)} \right)\\
\\
& &\qquad \displaystyle \times \left |  \frac{  (-1+a)a^{n_1+n_2}(a^{n_1}-a^{n_2})   }{\left( x_1a^{n_2}(-1+a^{n_1})-x_2 a^{n_1}(-1+a^{n_2})\right)^2}  \right|\mathrm{d}a.
\end{array}
\end{equation}

\section{PDF of the steady state}\label{sec_steady_state}

An important  issue in dealing with Beverton-Holt model is to determine the steady state.
From the deterministic theory it is known that its steady state is given by \eqref{estacionario}. As $\mathbb{P}\left[ \{ \omega \in \Omega:\, A(\omega)=1 \} \right]=0$, for all event $\omega\in \Omega$, the steady state of Beverton-Holt model  \eqref{Pielourand} is given by 
\begin{equation}\label{estacionariorand}
X_\infty=
\frac{A-1}{B}.
\end{equation}
This section is devoted to compute the PDF of \eqref{estacionariorand}. To achieve this objective, we define the mapping $\mathbf{r}$ by
\[
\begin{array}{ccccl}
y_1 & = & r_1\left( c,a,b \right) &=& \displaystyle \frac{a-1}{b},\\[0.2cm]
y_2 & = & r_2\left( c,a,b \right) &=& c,\\[0.2cm]
y_3 & = & r_3\left( c,a,b \right) &=& b,
\end{array}
\]
and using RVT technique (Theorem \ref{R.V.T.method2}), it is easy to check that the PDF corresponding  to $f_{Y_1,Y_2,Y_3}(y_1,y_2,y_3)$ is given by
\[
 f_{Y_1,Y_2,Y_3}(y_1,y_2,y_3) = f_{C,A,B}(y_2,y_1 y_3 + 1,y_3) \left| y_3 \right|.
\]
Finally, the PDF of the steady state \eqref{estacionariorand} is obtained marginalizing last expression, obtaining
\begin{equation}
f_{X_{\infty}}(x)=
\iint_{\mathbb{R}^2}
f_{C,A,B}
\left( c, x b+1,b  \right)
\left| b \right| \dd c\, \dd b.
\end{equation}

We will check via the example exhibited in next section that $f_1^{X}(x,n)$ given by \eqref{f1xn} converges to $f_{X_{\infty}}(x)$ as $n$ increases.

\section{Numerical experiments} \label{sec_ejemplos}

In this section we illustrate our theoretical results via several numerical examples. Computations have been carried out using the software Mathematica\textsuperscript{\tiny\textregistered}. 

\begin{figure}[htp]
\begin{center}
\includegraphics[scale=0.9]{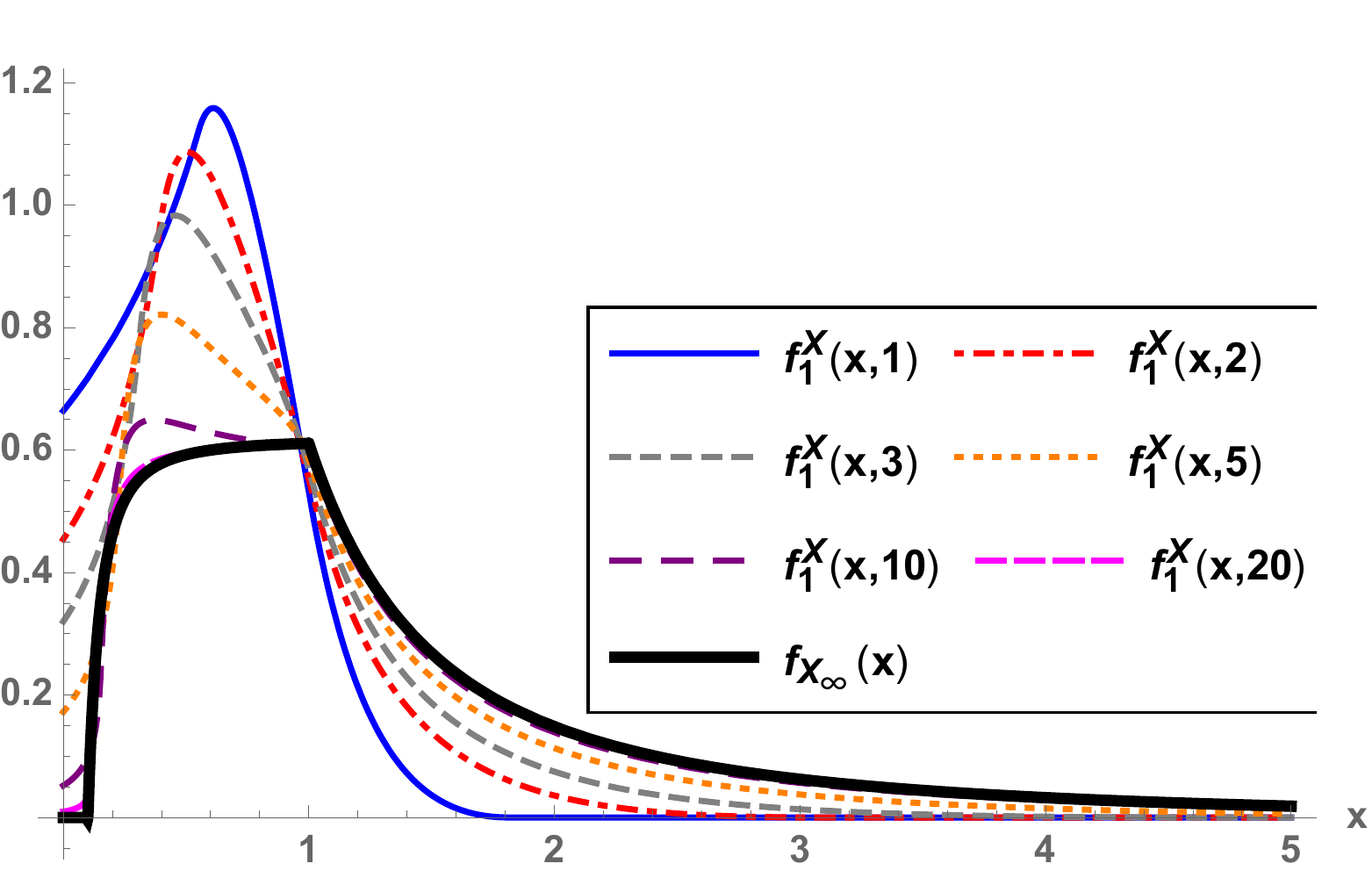}
\end{center}\caption{1-PDF, $f_1^X(x;n)$, of the solution, $X_n$, with $n\in\{1,2,3,5,10,20\}$ and the PDF, $f_{X_{\infty}}(x)$, of the equilibrium or steady state  $X_{\infty}$. Example \ref{ex1}.}
\label{pdf}
\end{figure}

\begin{figure}[htp]
\begin{center}
\includegraphics[scale=0.9]{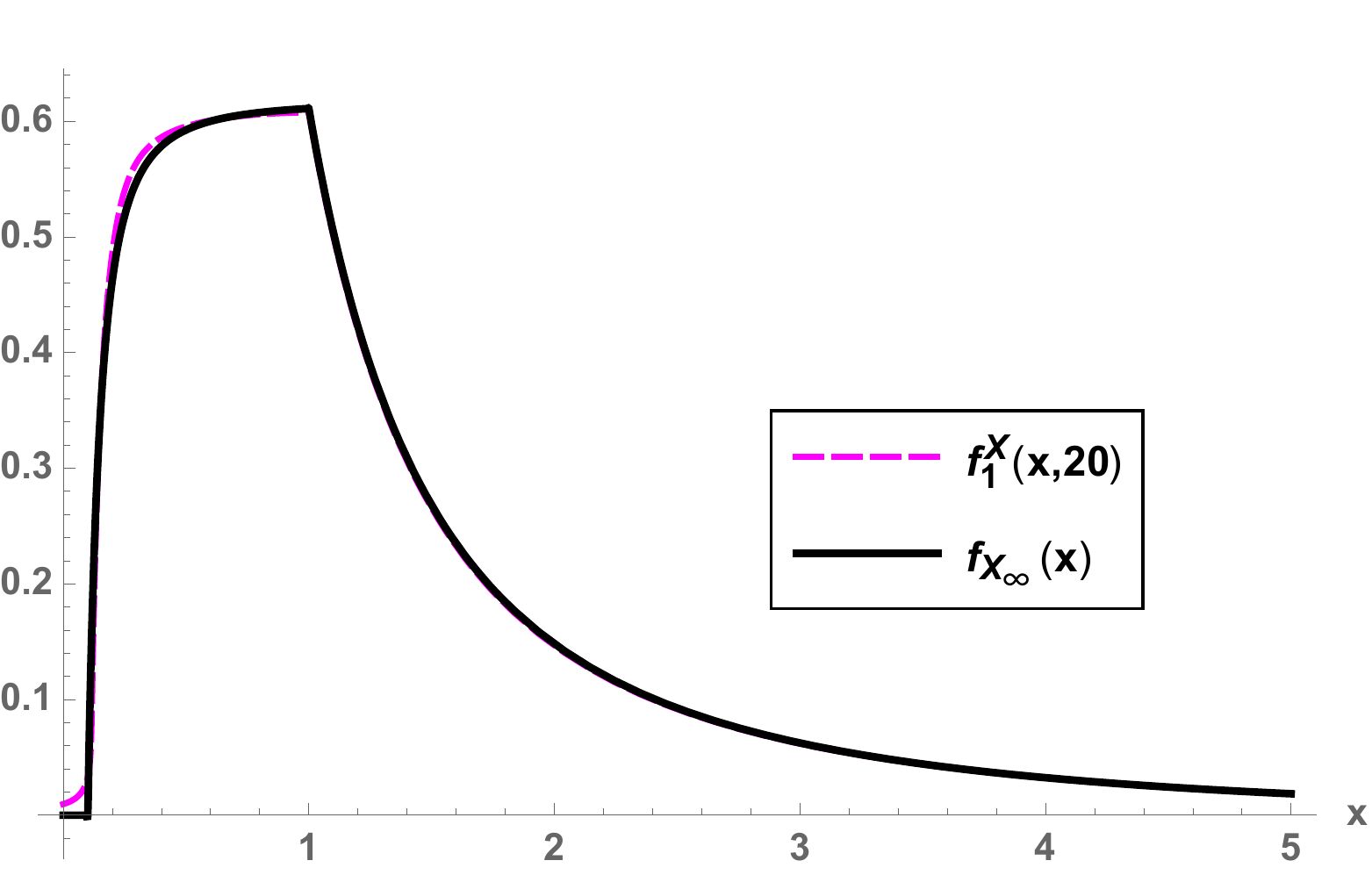}
\end{center}\caption{1-PDF, $f_1^X(x;n)$, of the solution, $X_{n}$ with $n=20$, and PDF, $f_{X_{\infty}}(x)$, of steady state $X_{\infty}$. Example \ref{ex1}.}
\label{pdfe}
\end{figure}

\begin{figure}[htp]
\begin{center}
\includegraphics[scale=0.55]{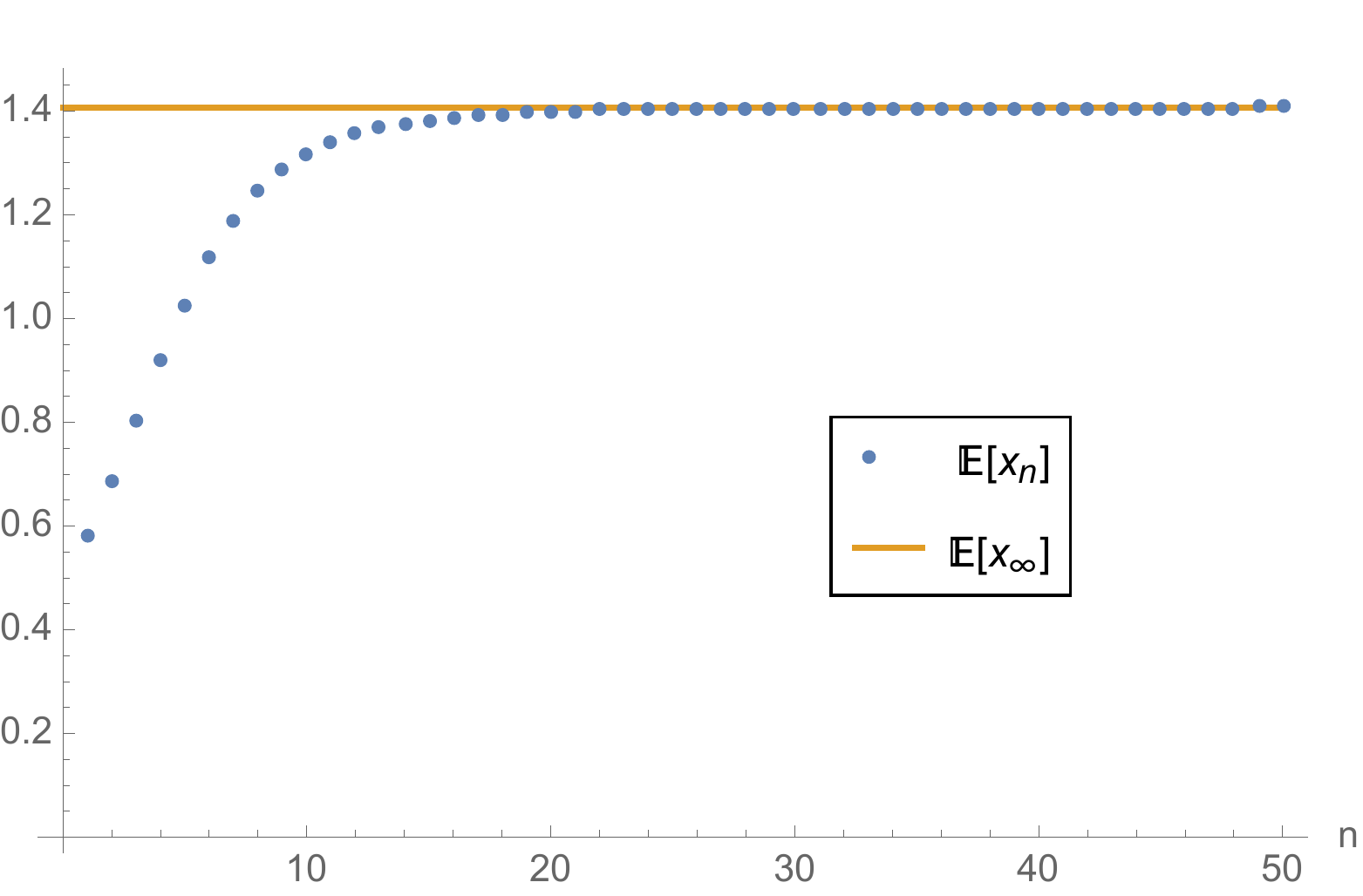}
\includegraphics[scale=0.55]{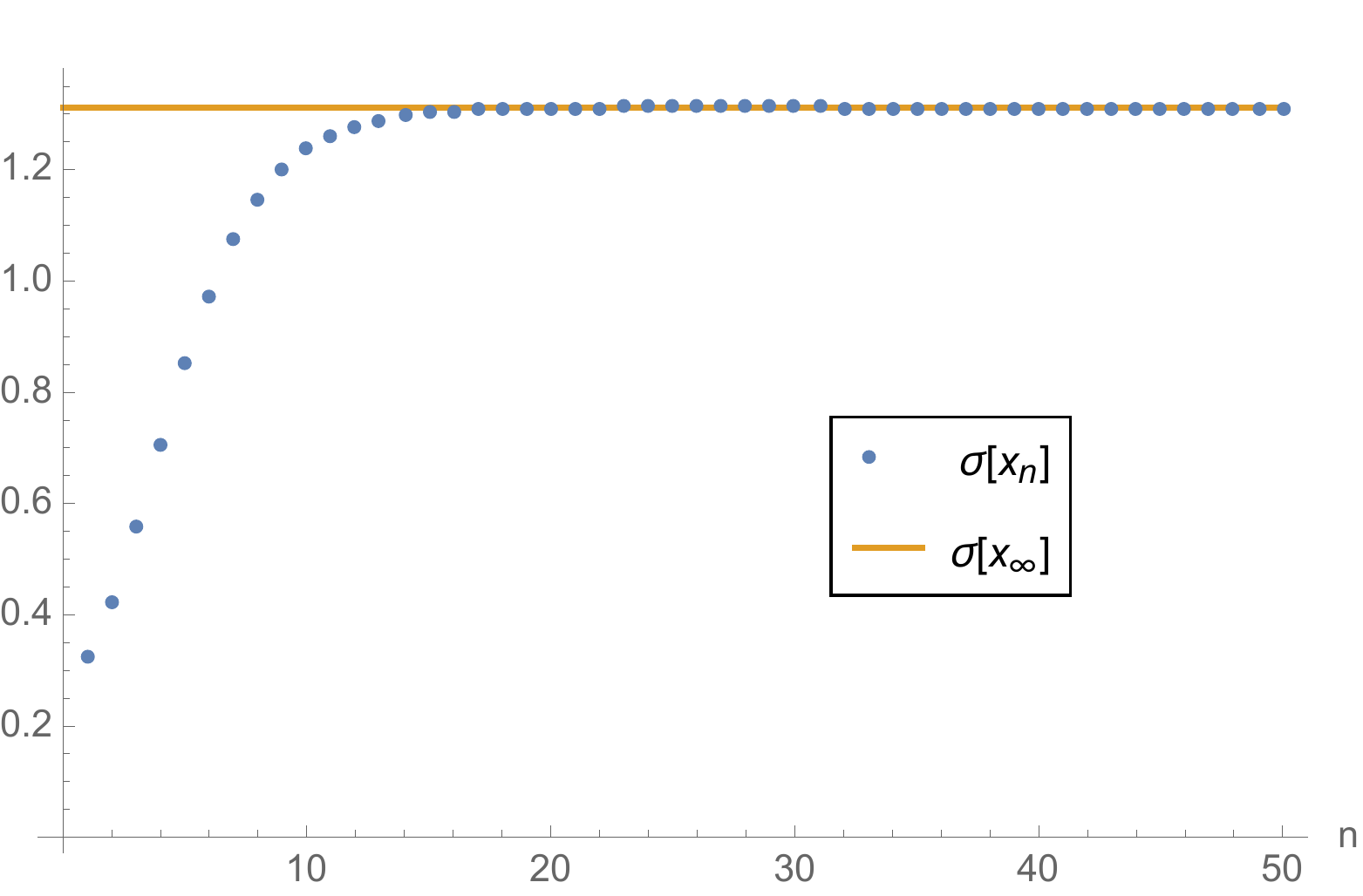}
\end{center}\caption{Left: Points represent the mean of $X_n$, $\mathbb{E}\left[X_n\right]$, for different $n\in \{1,2,\dots,50\}$. Solid line represents the mean,  $\mathbb{E}\left[X_{\infty}\right]$, of the steady state $X_{\infty}$. Right: Points represent the standard deviation of $X_n$, $\sigma\left[X_n\right]$, for different $n\in \{1,2,\dots,50\}$. Solid line represents the standard deviation, $\sigma\left[X_{\infty}\right]$, of the steady state $X_{\infty}$. Example \ref{ex1}.}
\label{media_1}
\end{figure}

\begin{figure}[htp]
\begin{center}
\includegraphics[scale=0.65]{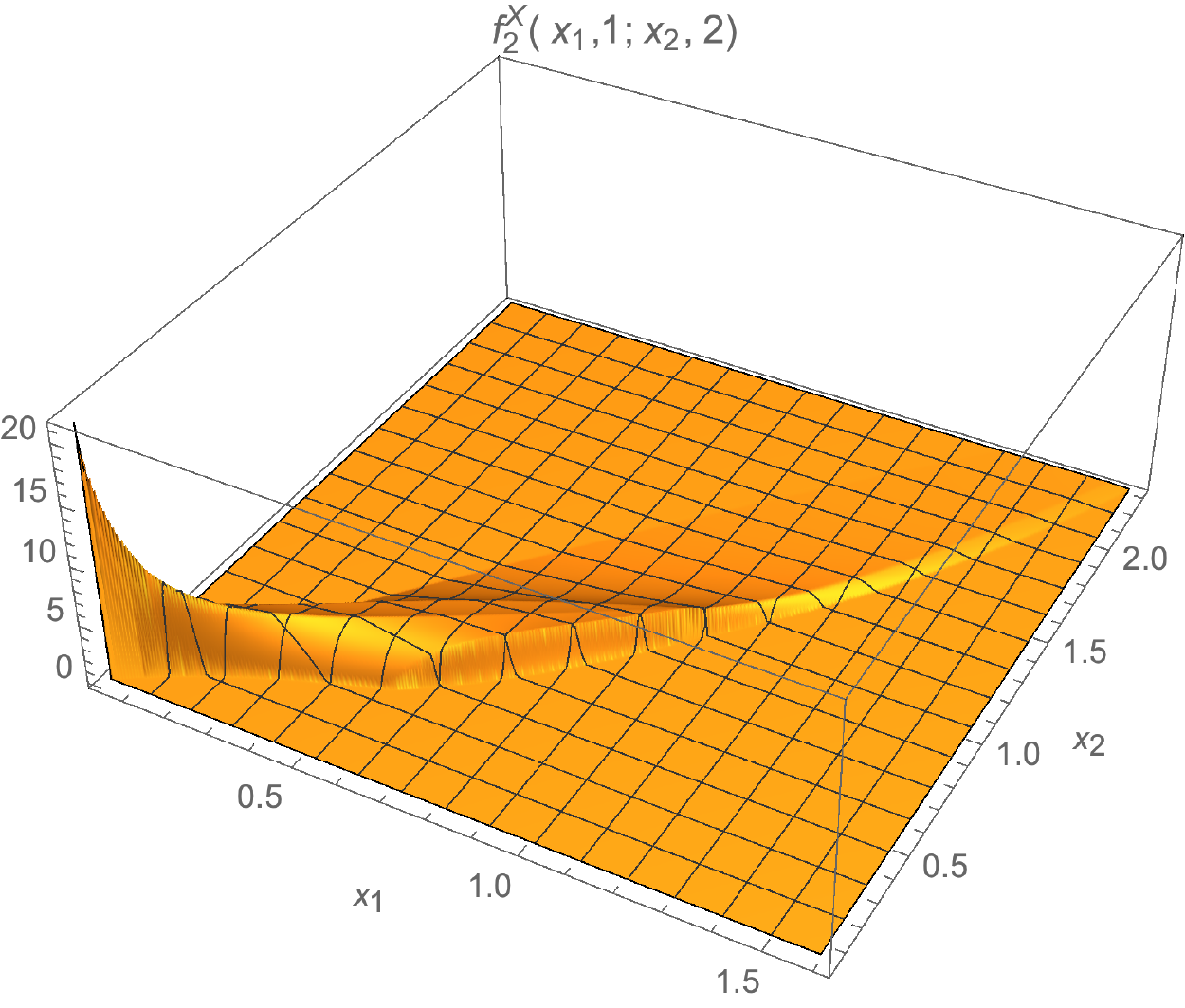}
\includegraphics[scale=0.65]{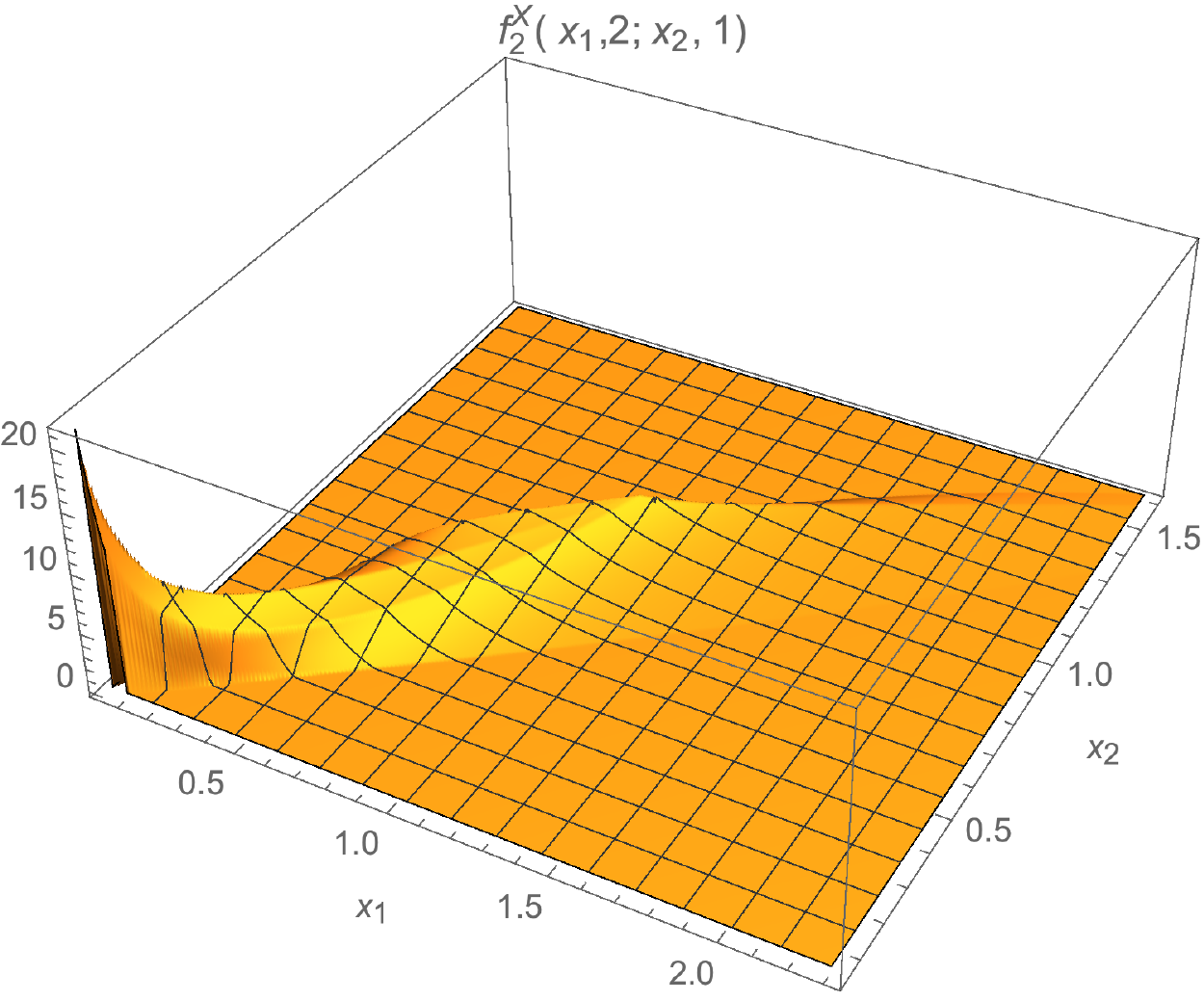}
\end{center}\caption{2-PDF, $f_2^X(x_1,n_1;x_2,n_2)$, of the solution, $X_n$. Left:  $n_1=1$ and $n_2=2$. Right:  $n_1=2$ and $n_2=1$. Example \ref{ex1}.}
\label{2PDF}
\end{figure}

\begin{figure}[htp]
\begin{center}
\includegraphics[scale=0.9]{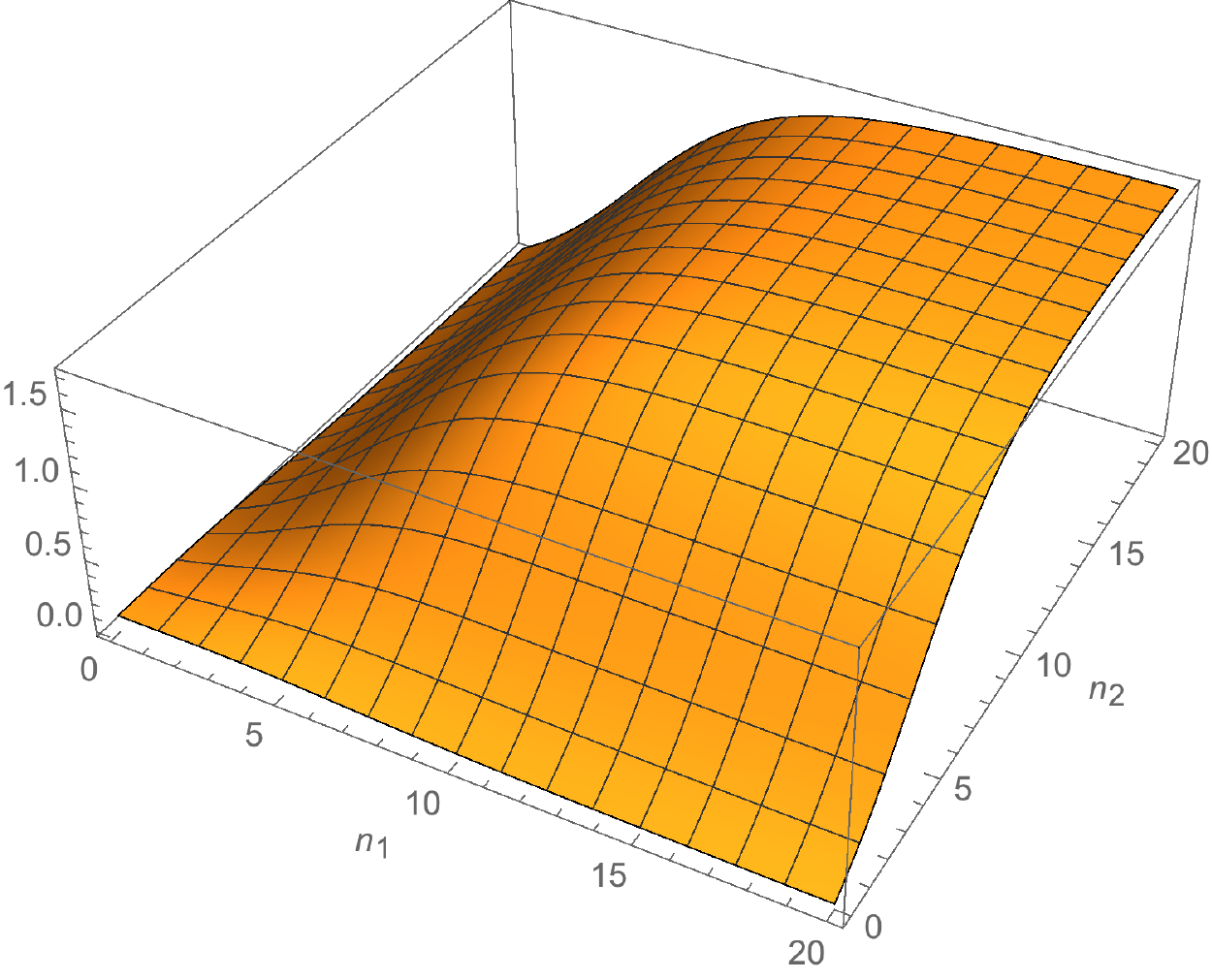}
\end{center}\caption{Covariance function, $C_X(n_1,n_2)$, of the solution, $X_n$, for the values of $(n_1, n_2)\in[0,20] \times [0,20]$. Example \ref{ex1}.}
\label{cov}
\end{figure}

\begin{example} \label{ex1} \normalfont

This section is addressed to illustrate our previous theoretical results  through an example.
With this aim, we need to choose the distributions of random variables $C$, $A$ and $B$. We will assume that these random variables are independent. Therefore, its joint PDF can be factorized as the product of each individual PDF,  $f_{C,A,B}(c,a,b)=f_{C}(c)f_{A}(a)f_{B}(b)$. As $X$ must 
be positive we have chosen for $C$ a uniform distribution in the interval $[0,1]$. 
For $A$ we have taken a  uniform distribution in the interval $[1.1,2]$ and for $B$ a  uniform distribution in the interval $[0.1,1]$

In Figure \ref{pdf}, the 1-PDF, $f_1^{X}(x,n)$, of the solution stochastic process, $X_n$,  for different values of $n\in \{1,2,3,5,10,20\}$ and the PDF, $f_{X_{\infty}}(x)$, of the equilibrium point, $X_{\infty}$,  are shown. As it has been previously pointed out, in this graphical representation one can observe that $f_1^{X}(x,n)$ tends to  $f_{X_{\infty}}(x)$ when $n$ increases. In Figure \ref{pdfe}, we have plotted separately the PDFs $f_1^{X}(x;20)$ and $f_{X_{\infty}}(x)$ to better highlight  this behaviour. In this figure we can observe that both PDFs match. 

In the left part of Figure \ref{media_1}, the mean of the solution stochastic process, $X_n$, and the threshold computed by the mean of the equilibrium have been represented. As it occurs with the 1-PDF, now we can observe that $\mathbb{E}\left[X_n\right]$ tends to $\mathbb{E}\left[X_{\infty}\right]$ as $n$ increases. In the right part a similar plot is presented for the standard deviation.

\begin{center}
\textbf{Figure 1}. 1-PDF, $f_1^X(x;n)$, of the solution, $X_n$, with $n\in\{1,2,3,5,10,20\}$ and the PDF, $f_{X_{\infty}}(x)$, of the equilibrium or steady state  $X_{\infty}$. Example \ref{ex1}.
\end{center}

\begin{center}
\textbf{Figure 2}. 1-PDF, $f_1^X(x;n)$, of the solution, $X_{n}$ with $n=20$, and PDF, $f_{X_{\infty}}(x)$, of steady state $X_{\infty}$. Example \ref{ex1}.
\end{center}

\begin{center}
\textbf{Figure 3}. Left: Points represent the mean of $X_n$, $\mathbb{E}\left[X_n\right]$, for different $n\in \{1,2,\dots,50\}$. Solid line represents the mean,  $\mathbb{E}\left[X_{\infty}\right]$, of the steady state $X_{\infty}$. Right: Points represent the standard deviation of $X_n$, $\sigma\left[X_n\right]$, for different $n\in \{1,2,\dots,50\}$. Solid line represents the standard deviation, $\sigma\left[X_{\infty}\right]$, of the steady state $X_{\infty}$. Example \ref{ex1}.
\end{center}

In Figure \ref{2PDF}, the 2-PDF of the solution stochastic process $X_n$ has been plotted in two cases: (1) $(n_1,n_2)=(1,2)$ and  $(n_1,n_2)=(2,1)$. Finally, in Figure \ref{cov}, the covariance surface, $C_{X}(n_1,n_2)$ of the solution stochastic process $X_n$ has been represented. This important deterministic function has been computed taking advantage of the explicit expressions of the 2-PDF $f_2^X(x_1,n_1;x_2,n_2)$ given in  \eqref{2PDFconj} together with expressions \eqref{correlacion}--\eqref{covarianza}.

\begin{center}
\textbf{Figure 4}. 2-PDF, $f_2^X(x_1,n_1;x_2,n_2)$, of the solution, $X_n$. Left:  $n_1=1$ and $n_2=2$. Right:  $n_1=2$ and $n_2=1$. Example \ref{ex1}.
\end{center}

\begin{center}
\textbf{Figure 5}. Covariance function, $C_X(n_1,n_2)$, of the solution, $X_n$, for the values of $(n_1, n_2)\in[0,20] \times [0,20]$. Example \ref{ex1}.
\end{center}

\end{example}

\begin{figure}[htp]
\begin{center}
\includegraphics[scale=0.9]{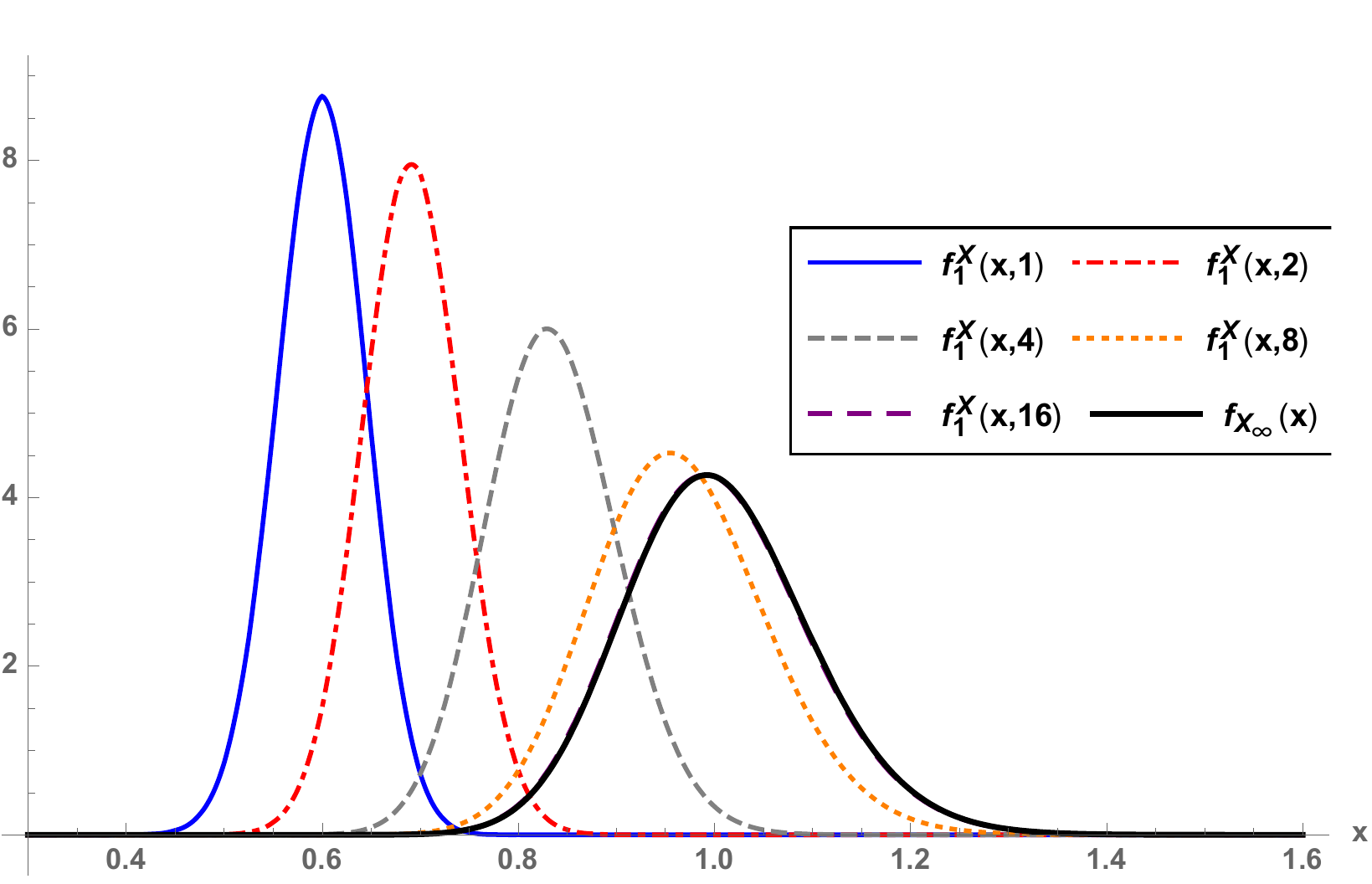}
\end{center}\caption{1-PDF, $f_1^X(x,n)$, of the solution, $X_n$, with $n\in\{1,2,4,8,16\}$ and the PDF, $f_{X_{\infty}}(x)$, of the equilibrium or steady state  $X_{\infty}$. Example \ref{ex2}.}
\label{pdf2}
\end{figure}

\begin{figure}[htp]
\begin{center}
\includegraphics[scale=0.9]{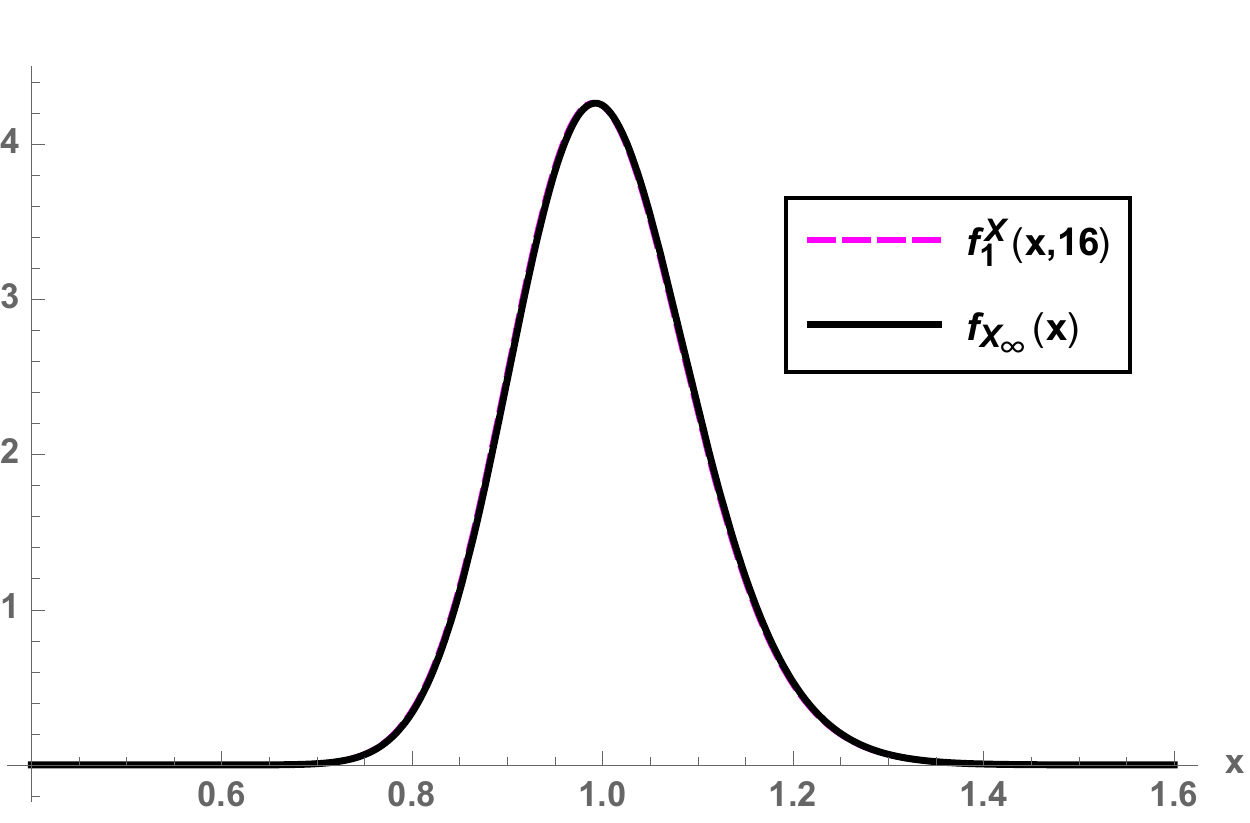}
\end{center}\caption{1-PDF, $f_1^X(x,n)$, of the solution, $X_{n}$, with $n=16$, and PDF, $f_{X_{\infty}}(x)$, of the equilibrium or steady state  $X_{\infty}$. Example \ref{ex2}.}
\label{pdfe2}
\end{figure}

\begin{figure}[htp]
\begin{center}
\includegraphics[scale=0.55]{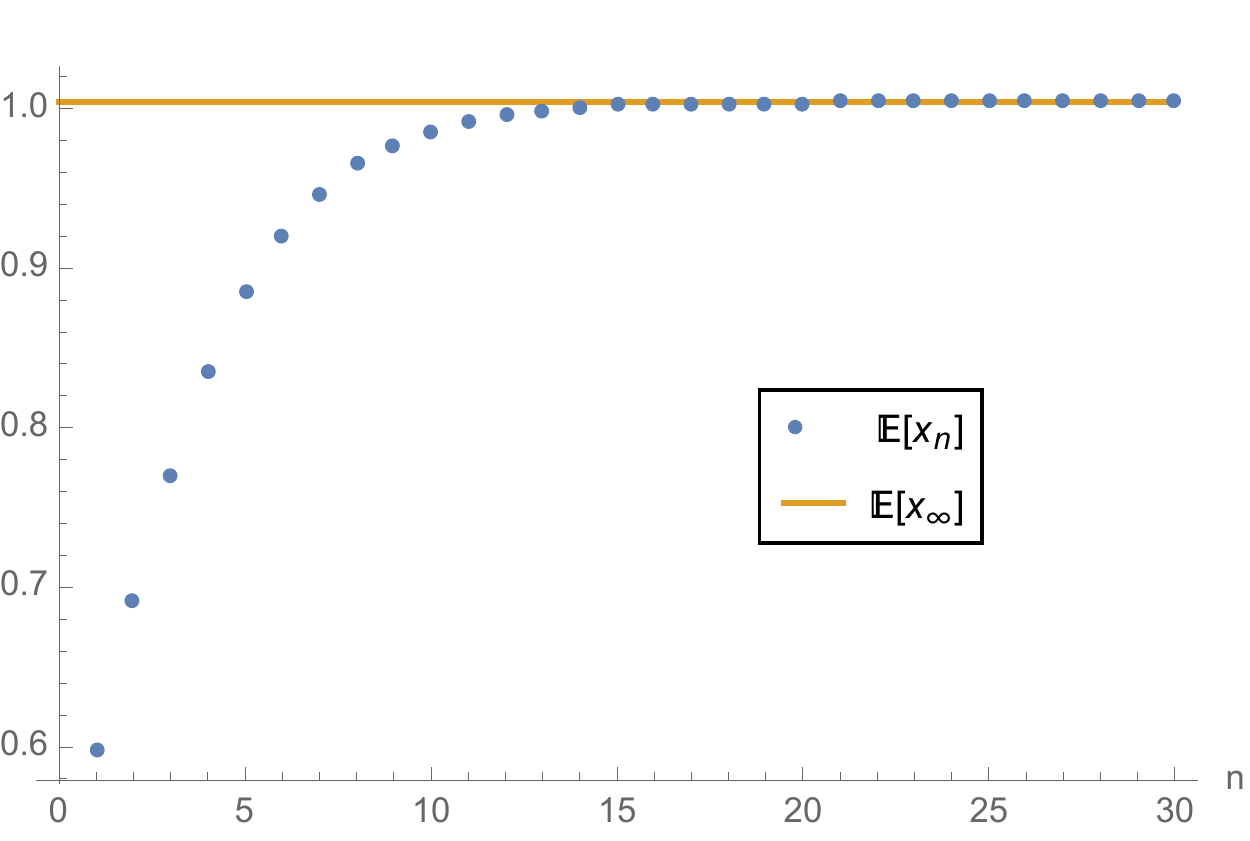}
\includegraphics[scale=0.55]{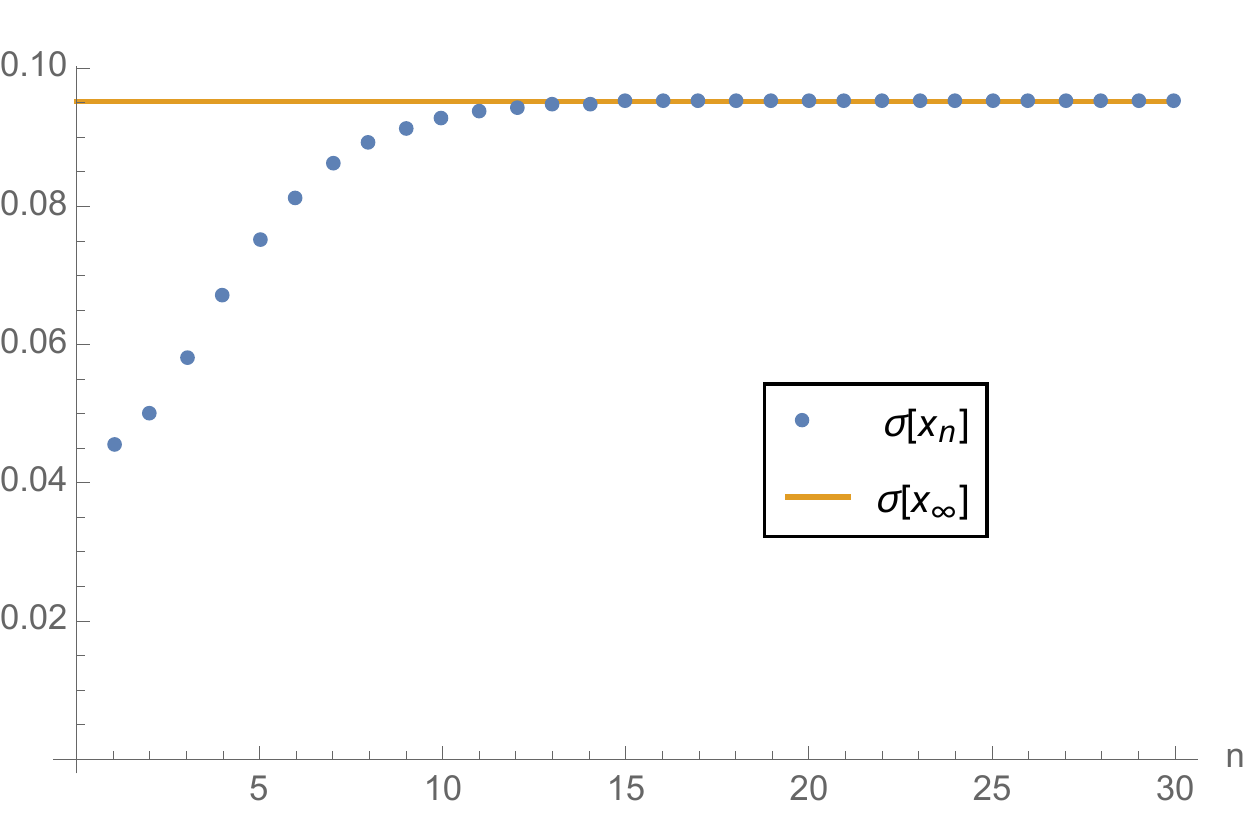}
\end{center}\caption{Left: Points represent the mean, $\mathbb{E}\left[X_n\right]$, of $X_n$ for  $n\in \{1,2,\dots,30\}$. Solid line  represents the mean, $\mathbb{E}\left[X_{\infty}\right]$, of the steady state, $X_{\infty}$. Right: Points represent the standard deviation, $\sigma\left[X_n\right]$, of $X_n$ for   $n\in \{1,2,\dots,30\}$. Solid line represents the standard deviation, $\sigma\left[X_{\infty}\right]$, of the steady state $X_{\infty}$. Example \ref{ex2}.}
\label{media_2}
\end{figure}

\begin{figure}[htp]
\begin{center}
\includegraphics[scale=0.5]{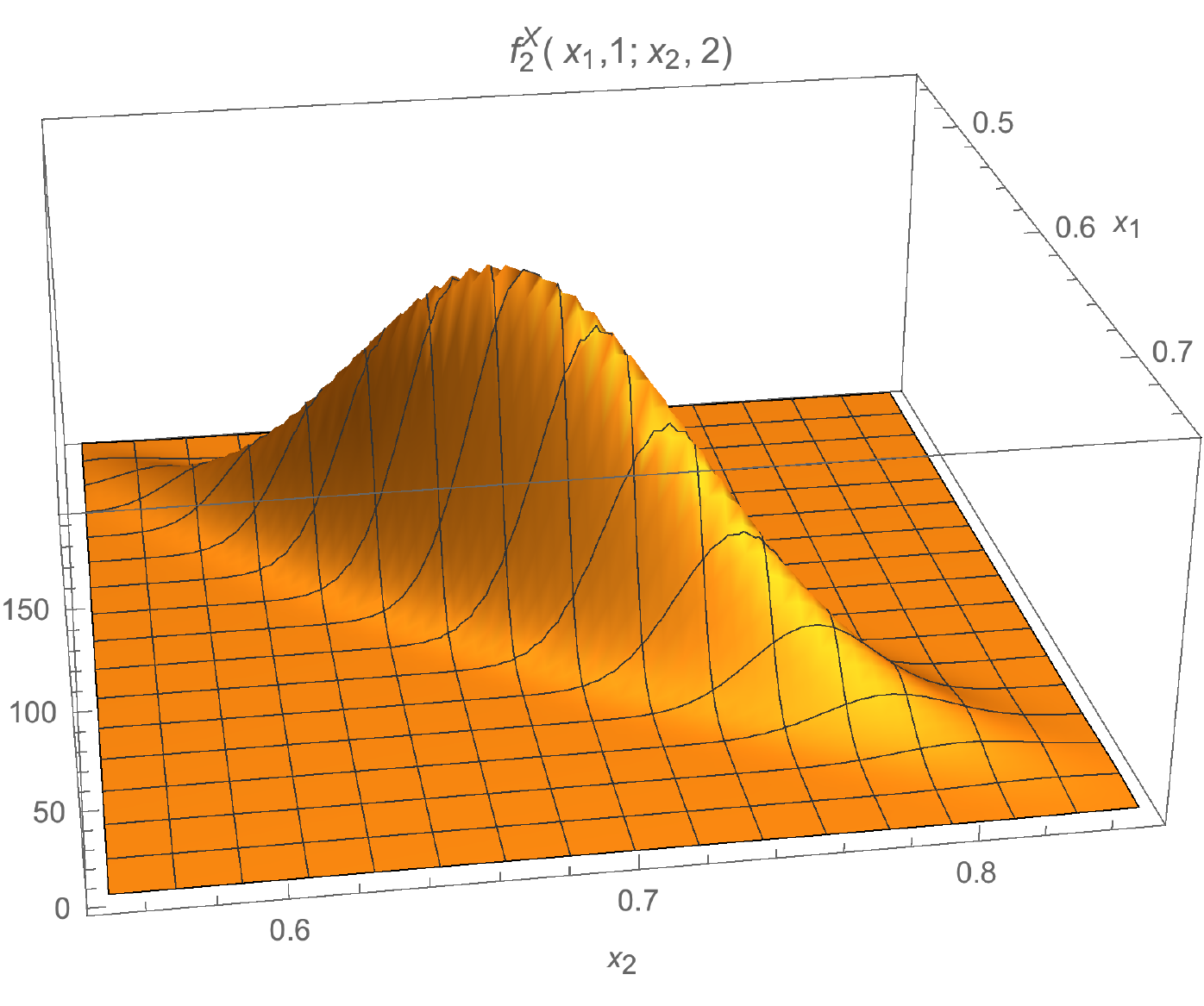}
\includegraphics[scale=0.5]{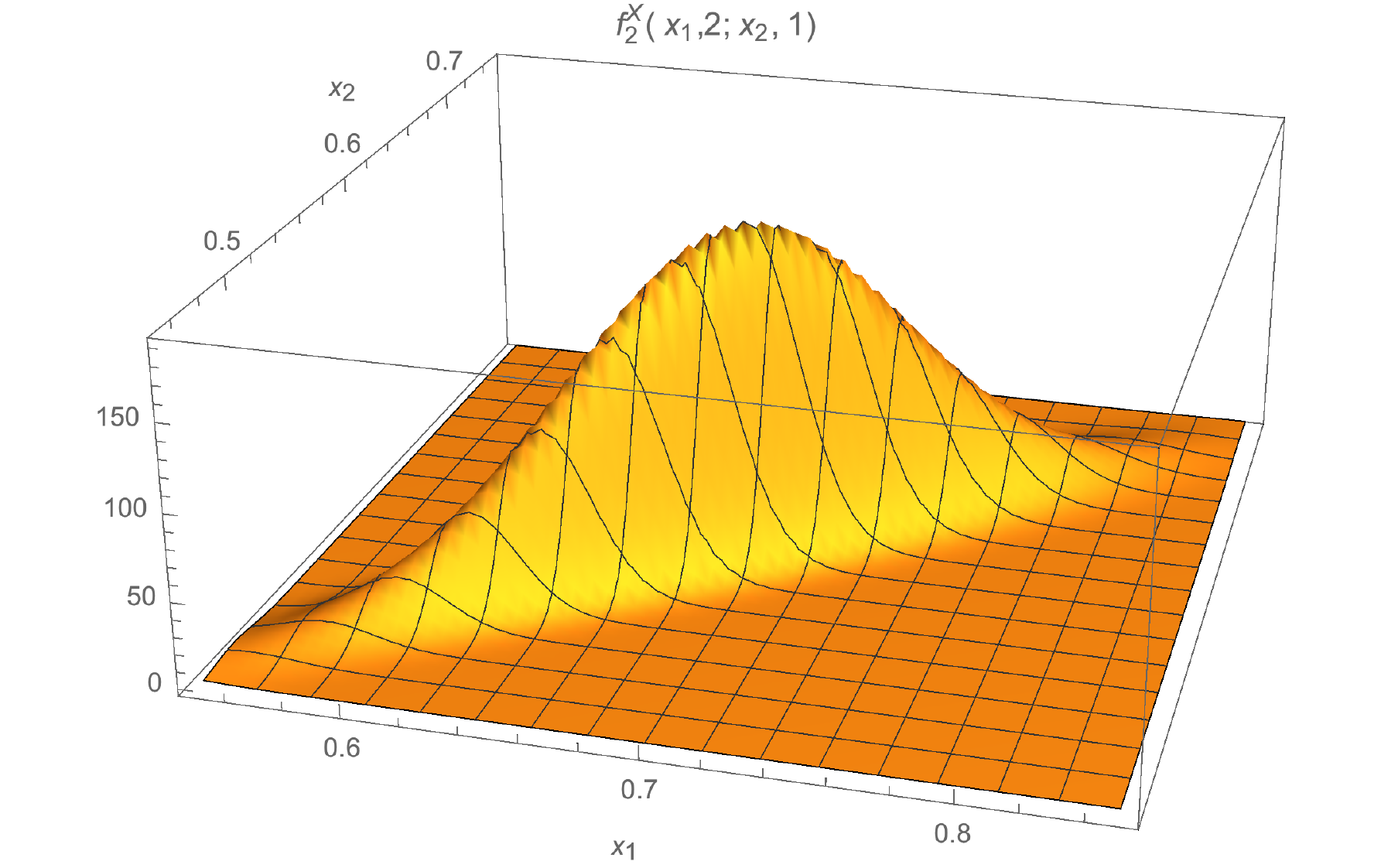}
\end{center}\caption{2-PDF, $f_2^X(x_1,n_1;x_2,n_2)$, of the solution, $X_n$. Left:  $n_1=1$ and $n_2=2$. Right: $n_1=2$ and $n_2=1$. Example \ref{ex2}.}
\label{2PDF2}
\end{figure}

\begin{figure}[htp]
\begin{center}
\includegraphics[scale=0.9]{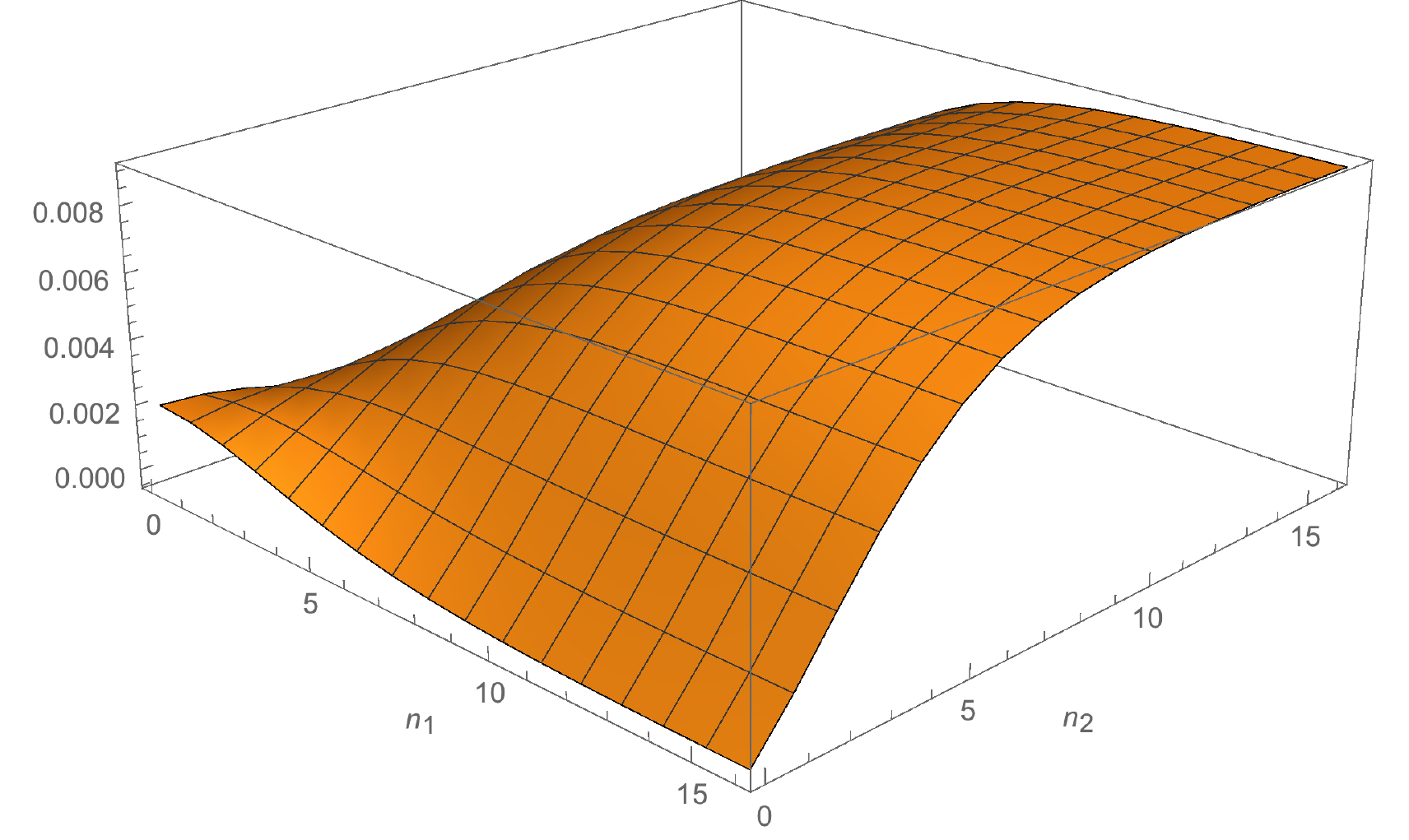}
\end{center}\caption{Covariance function, $C_X(n_1,n_2)$, of the solution, $X_n$, for the values of $(n_1, n_2)\in[0,16] \times [0,16]$. Example \ref{ex2}.}
\label{cov2}
\end{figure}

\begin{example} \label{ex2} \normalfont
In contrast to Example \ref{ex1}, where we have considered that the involved input data  are statistically independent, in this example we will illustrate our theoretical findings  assuming that  $A$, $B$ and $C$ are dependent random variables.  Specifically, hereinafter we will assume that the random vector $(C,A,B)$  has  a truncated Gaussian distribution $(C,A,B)\sim \text{N}_T(\mu,\Sigma)$, where $T=[0,1]\times [1.1,2 ]\times[0,1]$ is the domain of truncation, and the parameters $\mu$ and $\Sigma$, that represent the mean vector and the variance-covariance matrix, respectively, are  given by
\[
\mu=\left[ \begin{array}{c}
0.5\\ 1.5 \\ 0.5
\end{array}\right],\quad \Sigma=\frac{1}{500}\left[\begin{array}{ccc}
1 &0.1&0.2\\
0.1&0.9& 0.3\\
0.2& 0.3 &0.8
\end{array}\right].
\]

In Figure \ref{pdf2}, the 1-PDF, $f_1^{X}(x,n)$, of the solution stochastic process, $X_n$,  for different values of $n\in \{1,2,4,8,16\}$ and the PDF, $f_{X_{\infty}}(x)$, of the equilibrium point, $X_{\infty}$,  are shown. As it has been previously pointed out, in this graphical representation one can observe that $f_1^{X}(x,n)$ tends to  $f_{X_{\infty}}(x)$ as $n$ increases. In Figure \ref{pdfe2}, we have plotted separately the PDFs $f_1^{X}(x,16)$ and $f_{X_{\infty}}(x)$ to show better this behaviour. In this figure we can observe that both PDFs match. 

In Figure \ref{media_2} (left),  both the mean of the solution stochastic process, $X_n$, and the threshold computed by the mean of the equilibrium have been represented. As it occurs with the 1-PDF, now we can observe that $\mathbb{E}\left[X_n\right]$ tends to $\mathbb{E}\left[X_{\infty}\right]$ as $n$ increases. In Figure \ref{media_2} (right), a similar plot is shown for the standard deviation.

\begin{center}
\textbf{Figure 6}. 1-PDF, $f_1^X(x,n)$, of the solution, $X_n$, with $n\in\{1,2,4,8,16\}$ and the PDF, $f_{X_{\infty}}(x)$, of the equilibrium or steady state  $X_{\infty}$. Example \ref{ex2}.
\end{center}

\begin{center}
\textbf{Figure 7}. 1-PDF, $f_1^X(x,n)$, of the solution, $X_{n}$, with $n=16$, and PDF, $f_{X_{\infty}}(x)$, of the equilibrium or steady state  $X_{\infty}$. Example \ref{ex2}.
\end{center}

\begin{center}
\textbf{Figure 8}. Left: Points represent the mean, $\mathbb{E}\left[X_n\right]$, of $X_n$ for  $n\in \{1,2,\dots,30\}$. Solid line  represents the mean, $\mathbb{E}\left[X_{\infty}\right]$, of the steady state, $X_{\infty}$. Right: Points represent the standard deviation, $\sigma\left[X_n\right]$, of $X_n$ for   $n\in \{1,2,\dots,30\}$. Solid line represents the standard deviation, $\sigma\left[X_{\infty}\right]$, of the steady state $X_{\infty}$. Example \ref{ex2}.
\end{center}

In Figure \ref{2PDF2}, the 2-PDF of the solution stochastic process, $X_n$, has been plotted in two cases: (1) $(n_1,n_2)=(1,2)$ and  $(n_1,n_2)=(2,1)$. Finally, in Figure \eqref{cov2}, the covariance surface, $C_{X}(n_1,n_2)$, of the solution stochastic process, $X_n$, has been represented on the square $(n_1, n_2)\in[0,16] \times [0,16]$. This important deterministic function has been computed taking advantage of the explicit expression of the 2-PDF, $f_2^X(x_1,n_1;x_2,n_2)$, given in  \eqref{2PDFconj} together with expressions \eqref{correlacion}--\eqref{covarianza}.

\begin{center}
\textbf{Figure 9}. 2-PDF, $f_2^X(x_1,n_1;x_2,n_2)$, of the solution, $X_n$. Left:  $n_1=1$ and $n_2=2$. Right: $n_1=2$ and $n_2=1$. Example \ref{ex2}.
\end{center}

\begin{center}
\textbf{Figure 10}. Covariance function, $C_X(n_1,n_2)$, of the solution, $X_n$, for the values of $(n_1, n_2)\in[0,16] \times [0,16]$. Example \ref{ex2}.
\end{center}

\end{example}

\section{Conclusions} \label{sec_Conclusiones}

In  this  paper we have randomized the Beverton-Holt model. Then, we have provided a full probabilistic description of its solution stochastic process  under very general assumptions on random input data. That description has been made through the 1-PDF and the 2-PDF of the discrete solution stochastic process. A full probabilistic description through its PDF is also given to the steady state. Finally, some numerical examples illustrating our theoretical results have been shown.

\section*{Acknowledgements}
This work has been partially supported by the Ministerio de Econom\'{i}a y Competitividad grant MTM2017-89664-P. Ana Navarro Quiles acknowledges the postdoctoral contract financed by DyCon project funding from the European Research Council (ERC) under the European Union’s Horizon 2020 research and innovation programme (grant agreement No 694126-DYCON).

\section*{Conflict of Interest Statement}
The authors declare that there is no conflict of interests regarding the publication of this article.




\newpage

\newpage

\newpage

\newpage

\newpage

\newpage

\newpage

\newpage

\newpage

\newpage

\begin{thebibliography}{a}

\bibitem{molles2007ecology}
Molles, M. C. and Cahill, J. F. \textit{Ecology: Concepts and Applications}. McGraw-Hill Ryerson,  2007.

\bibitem{ver}
Verhulst, P. F. \textit{Recherches math\'ematiques sur la loi d'accroissement de la population}. Nouvelles m\'em. de l'Academie Royale des Sciences et Belles-Lettres de Bruxelles, 18, 1--41 (1845). 

\bibitem{ver2}
Verhulst, P. F. \textit{Deuxi\`eme m\'emoire sur la loi d'accroissement de la population}. Nouvelles m\'em. de l'Academie Royale des Sciences et Belles-Lettres de Bruxelles, 20, 1--32 (1845). 

\bibitem{KWASNICKI20135}
Kwasnicki W. \textit{Logistic growth of the global economy and competitiveness of nations}. Technological Forecasting and Social Change, 80(1) 50--76 (2013). Doi: 10.1016/j.techfore.2012.07.007.

\bibitem{Banasiak}
Banasiak J. \textit{Mathematical Modelling in One Dimension}. Cambrigde, 2013. 

\bibitem{Beverton}
Beverton, R.J.H.,  Holt, S.J. \textit{On the dynamics of exploited fish populations}. Fishery Investigations, Vol. 19. (Great Britain, Ministry of Agriculture, Fisheries, and Food). London: H.M. Stationery off.,1957.

\bibitem{De_La_Sen}
De la Sen M. \textit{The generalized Beverton-Holt equation and the control of populations}. Applied Mathematical Modelling, 32(11), 190--208  (2008). Doi: 10.1016/j.apm.2007.09.007.

\bibitem{CORTES2018190}
Cort\'es, J. C., Navarro-Quiles A.,  Romero J. V.,  Rosell\'o M. D. \textit{Computing the probability density function of non-autonomous first-order linear homogeneous differential equations with uncertainty}. Journal of Computational and Applied Mathematics, 337, 190--208  (2018). Doi: 10.1016/j.cam.2018.01.015.

\bibitem{CASABAN2017396}
Casab\'an, M. C., Cort\'es, J. C., Navarro-Quiles A.,  Romero J. V.,  Rosell\'o M. D., Villanueva, R.J. \textit{Computing probabilistic solutions of the Bernoulli random differential equation}. Journal of Computational and Applied Mathematics, 309, 396--407  (2017). Doi: 10.1016/j.cam.2016.02.034.

\bibitem{CORTES2017225}
Cort\'es, J. C., Navarro-Quiles A.,  Romero J. V.,  Rosell\'o M. D. \textit{Randomizing the parameters of a {Markov} chain to model the stroke disease: A technical generalization of established computational methodologies towards improving real applications}. Journal of Computational and Applied Mathematics, 324, 225--240  (2017). Doi: 10.1016/j.cam.2017.04.040.

\bibitem{CORTES2017150}
Cort\'es, J. C., Navarro-Quiles A.,  Romero J. V.,  Rosell\'o M. D. \textit{Full solution of random autonomous first-order linear systems of difference equations. Application to construct random phase portrait for planar systems}. Applied Mathematics Letters, 68, 150--156  (2017). Doi: 10.1016/j.aml.2016.12.015.

\bibitem{soong}
Soong T. T. \textit{Random Differential Equations in Science and Engineering}. Academic Press, New York, 1973.




\end{thebibliography}
\end{document}